\documentclass[12 pt]{amsart}

\usepackage{amsfonts,amsmath,amsthm,amssymb}
\usepackage{fullpage,lineno,setspace, nicefrac, subcaption, caption}
\usepackage{mathrsfs}
\usepackage{tikz}
\usetikzlibrary{automata,arrows,positioning,calc}
\usepackage[mathscr]{eucal}

\DeclareGraphicsExtensions{.png,.jpeg,.jpg,.pdf}
\graphicspath{{./Figures/}}

\newtheorem{theorem}{Theorem}[section]

\newcommand{\Def}{\overset{\textbf{def}}{=}}
\newcommand{\RR}{\mathbb{R}}

\newcommand{\abs}[1]{\lvert #1 \rvert}

\newcommand{\pp}{\partial}

\usepackage{amsopn}

\title{High-Order Numerical Method for 1D Non-local Diffusive Equation}
\author{Dawson Do}
\address{Department of Civil and Environmental Engineering\\
    University of California at Berkeley\\
    Berkeley, CA 94720}
\email{daws@berkeley.edu}
\author{Hossein Nick Zinat Matin}
\address{Department of Civil and Environmental Engineering\\
    University of California at Berkeley\\
    Berkeley, CA 94720}
\email{h-matin@berkeley.edu}
\author{Maria Laura Delle Monache}
\address{Department of Civil and Environmental Engineering\\
    University of California at Berkeley\\
    Berkeley, CA 94720}
\email{mldellemonache@berkeley.edu}
\begin{document}

\maketitle

\begin{abstract}
 In this paper we present a non-local numerical scheme based on the Local Discontinuous Galerkin method for a non-local diffusive partial differential equation with application to traffic flow.  In this model, the velocity is determined by both the average of the traffic density as well as the changes in the traffic density at a neighborhood of each point. We discuss nonphysical behaviors that can arise when including diffusion, and our measures to prevent them in our model. The numerical results suggest that this is an accurate method for solving this type of equation and that the model can capture desired traffic flow behavior. We show that computation of the non-local convolution results in $\mathcal{O}(n^2)$ complexity, but the increased computation time can be mitigated with high-order schemes like the one proposed.

  
\end{abstract}

\newenvironment{@abssec}[1]{%
     \if@twocolumn
       \section*{#1}%
     \else
       \vspace{.05in}\footnotesize
       \parindent .2in
         {\upshape\bfseries #1. }\ignorespaces 
     \fi}
     {\if@twocolumn\else\par\vspace{.1in}\fi}

\newenvironment{keywords}{\begin{@abssec}{Key words}}{\end{@abssec}}


\section{Introduction}
This paper focuses on the design of a non-local numerical scheme based on the Local Discontinuous Galerkin (LDG) method for non local diffusive model with application to traffic flow. 
Classical traffic flow models are based on the well-known  Lighthill-Whitham-Richards (LWR) model \cite{lighthill1955kinematic, richards1956shock}. It assumes that traffic can be modeled by a conservation law in which the mean traffic speed is a function of only the mean traffic density.  In the last decade, a new avenue of research has brought to light the possibility of adding non-local terms in the conservation laws to define physical characteristics, such as velocity,  on the average quantity within a ``neighborhood" of each point rather than locally. For example, in traffic flow, this is motivated by the drivers' ability to observe traffic status ahead of them and adjust their speed accordingly \cite{keimer2023nonlocal, blandin2016well, goatin2016well, bayen2022modeling, colombo2012class, huang2022stability}. 
Non-local terms have also been used in conservation law models in chemical engineering \cite{pflug2020emom}, in sedimentation \cite{betancourt2011nonlocal}, laser-cutting technology \cite{colombo2015nonlocal},  production network models \cite{gugat2016analysis} and population dynamics \cite{perthame2006transport}.

In addition to include non-local terms, transportation research has also focused on improving the LWR model by adding a diffusion term to the conservation law. This adjusts some nonphysical behavior of the LWR such as instantaneous speed change but, especially in the case of linear diffusion, violates other physical principles of traffic flow (i.e., the model can produce negative velocity) and was criticized by Daganzo \cite{daganzo1995requiem}. Later on, this was addressed with second order models  (\cite{aw2000resurrection,zhang2002non}) and in other papers focused on diffusive traffic flow models (see \cite{bonzani2000hydrodynamic,campos2021saturated,corli2021wavefronts,lighthill1955kinematic2}).

Several numerical methods have been developed to solve non-local conservation laws and advection-diffusion equations in general.
Authors have modified common schemes such as finite volume (FV) schemes (\cite{leveque2002finite}), high-order FV-WENO schemes (\cite{harten1997uniformly, shu1988efficient, shu1989efficient}), and discontinuous Galerkin (DG) schemes (\cite{cockburn2012discontinuous,cockburn1989tvb2,cockburn1989tvb}) to compute the non-local terms. Paper \cite{goatin2016well} studied non-local FV methods for traffic flow and the authors in \cite{gottlich2015discontinuous} compare FV and DG schemes for non-local material flow models in 2D. In \cite{chalons2018high}, the authors evaluate the performance of high-order FV-WENO schemes and DG schemes for non-local conservation laws in 1D, specifically for traffic flow and sedimentation models. Methods for solving advection-diffusion equations include implicit-explicit Runge-Kutta (RK) \cite{pareschi2005implicit} and the aforementioned FV-WENO. Additionally, DG methods for advection-diffusion equations originate from the LDG method, introduced by \cite{cockburn1998local}. To the best of our knowledge the only other instance of an LDG scheme for traffic flow applications is \cite{matin2023Nonlinear}. Notably, LDG is a natural higher-order generalization of the FV method, the scheme typically used in macroscopic traffic flow applications (see \cite{lebacque1996godunov}). This reason, along with LDG's stability for solving advection-dominated advection-diffusion equations, makes it a suitable choice for the proposed model. Furthermore, the computational results in \cite{ chalons2018high, gottlich2015discontinuous} encourage the use of high-order schemes, which are invaluable for efficiently simulating non-local models. 

In this paper, we focus on the design of a novel LDG method that is able to to handle diffusion and non-local terms. We do so,  by first introducing a new non-local diffusive traffic flow model that satisfies the physical principles of traffic flow. And then, developing a non-local LDG method that unifies the numerical scheme in \cite{chalons2018high} with LDG \cite{cockburn1998local}.


The paper is organized as follows. In Section \ref{sec:model}, we introduce the non-local diffusive model, discussing the physical motivation and considerations we have made in order to preserve realistic behavior. In Section \ref{sec:methods}, we derive the LDG-based scheme used to simulate the model. In Section \ref{sec:model_prop}, we validate the model and simulation by comparing our proposed model with conventional models. Finally, in Section \ref{sec:comp_results}, we present computational results. 

\section{Non-local diffusive traffic flow model}
\label{sec:model}
In this paper, we are concerned with non-local diffusive conservation laws with application in traffic flow dynamics, in the form of
\begin{equation}\label{E:main}
    \begin{cases}
        \rho_t + \partial_x [D(\rho) U(\hat \rho * \mathcal K_\gamma)] = 0 &, (t, x) \in \RR_+ \times \RR \\
        \rho(0, x) = \rho_\circ(x) 
    \end{cases}
\end{equation}
where,
\begin{equation*}
     D(\rho) \Def \rho(1 - \rho), \quad U(r) \Def 1 - r , \, \text{for any $\rho, r \in [0, 1]$}.
\end{equation*}
The velocity function $U(\cdot)$ can be any reasonable decreasing function of density. The function $(\rho, r) \mapsto (1 - \rho) U(r)$ represents the average velocity and $(\rho, r) \mapsto D(\rho)U(r)$ is the flux function. Furthermore, for a fixed $\kappa >0$, we define
\begin{equation} \label{E:fictitious_density}
    \hat \rho \Def \rho + \kappa \rho(1 - \rho) \Psi(\pp_x\rho), \quad \Psi(u) \Def \tanh(u), \, \forall u \in \RR.
\end{equation}
Let us recall that non-local models are proposed to improve the LWR model by considering the average density in determining the velocity rather than the local information. In the same spirit, the main idea for defining the model \eqref{E:main} and \eqref{E:fictitious_density} is that in reality, the density felt by the driver might be different from the traffic density $\rho$ and in particular dependent on $\pp_x \rho$; see \cite{bonzani2000hydrodynamic, de1999nonlinear}. 
Therefore, we calculate the function $U(\cdot)$, such that it considers the perceived density, i.e. $\hat \rho$. In particular, if $\pp_x \rho$ is positive, drivers feel a larger density as in \eqref{E:fictitious_density}, and similarly for the negative gradient of the density they feel less density than the average density $\rho$. On the other hand, we define the function $\Psi$ which is monotonically increasing, smooth and bounded, and hence the partial derivatives are prevented from growing unbounded. In particular, 
\begin{equation}\label{E:Psi_function}
    \Psi(\pp_x \rho) \to \pm 1 , \quad \text{as $\pp_x \rho \to \pm \infty$}
\end{equation}
From the traffic flow perspective, such boundedness ensures that the diffusion term is not the dominant term in defining the velocity. Let us next show that from the physical standpoint, $\hat \rho$ is a (perceived) density. 
\begin{theorem}
Let $\hat \rho$ be defined as in \eqref{E:fictitious_density} and $\kappa \in (0, 1)$. Then, $\hat \rho \in [0, 1]$.
\end{theorem}
\begin{proof}
First, it should be noted that when $\pp_x \rho$ increases (decreases) then $\hat \rho$ also increases (decreases) respectively. Equation \eqref{E:Psi_function} implies that as $\pp_x \rho \to - \infty$, the resulting quadratic equation $\rho -\kappa \rho(1 - \rho) - 1$ has the roots $\rho_1 = -\nicefrac{1}{\kappa}$ and $\rho_2 = 1$ and hence for $\rho \in (0, 1)$ and $\kappa \in (0, 1)$ the claim follows. 

For $\pp_x \rho \to \infty$ the quadratic equation $\rho +\kappa \rho(1 - \rho) - 1$ has the roots $\rho_1 = \nicefrac{1}{\kappa}$ and $\rho_2 = 1$ which means for $\kappa \in (0 ,1)$ the claim follows. In fact, the same argument directly can be used for $\abs{\pp_x \rho} < \infty$ by scaling $\kappa \in (0, 1)$ with respect to $\tanh(\pp_x \rho) \in (-1, 1)$. This completes the proof.
\end{proof}
Next, we elaborate on the role of the convolution term. In particular,
\begin{equation*}
 (\hat \rho * \mathcal K_\gamma)(t,x) \Def \int_{x}^{x + \gamma} \mathcal K_\gamma(y - x) \hat \rho(t, y) dy   
\end{equation*}
for a kernel $\mathcal K_\gamma$ presents a weighted average of densities $\hat \rho$ in a neighborhood of $x$. The length of this neighborhood is proportional to $\gamma$. To consider the weighted average of the surrounding density, the Kernel $\mathcal K_\gamma \in C_c^1([0, \gamma]; \RR_+)$ is defined to be decreasing so it vanishes outside the $\gamma$-neighborhood, and $\int_0 ^\gamma  \mathcal K_\gamma(x) dx= 1$. In this work, we consider the linear kernel 
\begin{equation} \mathcal K_\gamma (x) = \tfrac 2\gamma (1 - \tfrac x \gamma).
\end{equation}
Let
\begin{equation}\label{E:flux}
    q_F \Def D(\rho) U(\hat \rho * \mathcal K_\gamma).
\end{equation}
Then $q_F \ge 0$ (regardless of the growth rate of $\pp_x \rho$); i.e. by the construction of the problem the flow $q_F$ always moves in the positive direction even if $\pp_x \rho \to \pm \infty$ which from the application point of view this implies that the traffic moves in the correct direction. In addition, the model can be expanded to
\begin{equation*}
    \begin{split}
        \rho_t + \pp_x \rho D'(\rho) U(\hat \rho * \mathcal K_\gamma) = - D(\rho) \pp_x[\hat \rho * \mathcal K_\gamma] U'(\hat \rho * \mathcal K_\gamma)
    \end{split}
\end{equation*}
where $U'(\cdot) \equiv \partial_r U(\cdot)$. The term $\pp_x[\hat \rho * \mathcal K_\gamma]$ represents the diffusion term. This implies traffic dynamics are influenced both by the weighted average density and the weighted average of changes in the density. From the technical point of view, $U' \le 0$ and hence the diffusion coefficient on the right-hand side will be positive and we have a forward nonlinear parabolic equation. Finally, $D(\rho)$ on the right-hand side ensures the diffusion term degenerates properly at $\rho = 0$ and $\rho= 1$. In other words, when $\rho \to 0$ or $\rho \to 1$, the diffusion term vanishes. Consequently, the velocity is not determined by the diffusion in extreme cases and the model does not experience negative velocity. 

\section{Numerical Method}
\label{sec:methods}

An LDG scheme is chosen to simulate model \eqref{E:main}\eqref{E:fictitious_density} for its ease of incorporating an approximation of $\partial_x\rho$ into the solution. Furthermore, a high-order scheme is desired due to the accuracy and efficiency even at large mesh sizes. In this section, we derive the non-local LDG scheme for our model. Then we analyze the complexity of the method.

\subsection{Non-local LDG Scheme}
Starting with our equation: 
let us consider $q_F$ as in \eqref{E:flux}. Then, with a slight abuse of notation, we can write
\begin{equation}
        \pp_t\rho + \pp_x[q_F(\rho,\hat\rho*\mathcal{K}_\gamma)] = 0.
\end{equation}
In addition, for any $t \in \RR_+$, we define $R\Def \hat\rho(t,\cdot)*\mathcal{K}_\gamma(\cdot)$, and note that it is a function of $\rho$ and $\pp_x\rho$. Similar to LDG \cite{cockburn1998local}, we can rewrite and split into two equations by letting $\pp_x\rho=\sigma$:
\begin{equation}
    \begin{cases}
        \pp_t\rho + \pp_x[q_F(\rho,\sigma)] = 0 \\
        \pp_x\rho =\sigma .
    \end{cases}
\end{equation}
Let the function $(t,x) \mapsto u_h(t, x)$ represent our approximate solution, within the function space of discontinuous polynomials, $V_h$, over our domain partitioned into intervals of $[x_{k-1},x_{k}]$. $V_h$ has a set of basis functions, $\phi_i^k$, for each domain partition $k$, where $i\in\{0\dots p\}$ and $p$ is the highest polynomial degree of $V_h$. Basis functions are chosen such that they have the property that $\phi_i^k(x_i^k) = 1,\,\forall i$, where $x_i^k$ are basis points for each domain partition. In this paper, we use the Legendre basis polynomials defined on $p+1$ Chebyshev nodes as the basis points \cite{cockburn2001runge}. Notably, $x_{k} = x_p^{k}=x_0^{k+1}$. $u_h$ is defined as:
\begin{equation}
    u_h(x)=\sum_{k=1}^{n}u^k_h(x)=\sum_{k=1}^{n}\sum_{i=0}^p u_i^k\phi_i^k(x)
    \label{eq:approx_sol}
\end{equation}
where each $u_i^k$ is a real-valued coefficient. There is a similar definition for $\sigma_h$:
\begin{equation}
    \sigma_h(x)=\sum_{k=1}^{n}\sum_{i=0}^p \sigma_i^k\phi_i^k(x).
\end{equation}
Note that this means $u_h(x_i^k)=u_i^k$, and similarly for $\sigma_h$. Under the Galerkin formulation, we look for $u_h,\sigma_h\in V_h$ such that:
\begin{equation*}\int_\Omega \partial_t u_h vdx+\int_\Omega\partial_x (q_F(u_h,\sigma_h))vdx=0,\quad \forall v\in V_h\end{equation*} 
\begin{equation*}\int_\Omega \partial_t u_h\tau dx=\int_\Omega \sigma_h\tau dx,\quad \forall\tau\in V_h.\end{equation*}  
We can then set $v,\tau$ to each basis function of $V_h$, $\phi_i^k$, then integrate by parts. We have $\forall i \in\{0\dots p\},\,\forall k =\{1\dots n\}$:
\begin{equation*}
\int_{x_{k-1}}^{x_k} \partial_t u_h\phi_i^kdx+\Big[\hat q_F(x)\phi_i^k(x)\Big]_{x_{k-1}}^{x_k}-\int_{x_{k-1}}^{x_k}q_F(u_h,\sigma_h)\frac{d \phi_i^k}{dx}dx=0
\end{equation*} 
\begin{equation*}
\Big[\hat{u}(x)\phi_i^k(x)\Big]_{x_{k-1}}^{x_k}-\int_{x_{k-1}}^{x_k}u_h\frac{d \phi_i^k}{dx}dx=\int_{x_{k-1}}^{x_k} \sigma_h\phi_i^k dx
\end{equation*}
where $\hat{q}_F,\hat{u}$ are numerical fluxes at the boundary, since there is a discontinuity at each internal boundary. For this method, we use the Lax-Friedrich numerical flux for $\hat{q}_F$: 
\begin{equation*}
    \hat q_F(x_k) = Q(u_p^{k},u_0^{k+1},R^k) = \frac{1}{2}\left(\left(f\left(u_p^{k}\right) + f\left(u_0^{k+1}\right)\right)U\left(R^k\right)+\alpha\left(u_p^{k}-u_0^{k+1}\right)\right)
\end{equation*}
where $R^k=R(x_k)$ and $\alpha = \max\left|\pp_\rho(D(\rho)U(\hat \rho))\right|$. Consistent with the original LDG scheme, we use the numerical flux $\hat u(x_k)=u_0^{k+1}$. We then expand the internal boundary terms $\forall k,i$:
\begin{multline*}
    \int_{x_{k-1}}^{x_k} \partial_t u_h\phi_i^kdx-\int_{x_{k-1}}^{x_k}q_F(u_h,\sigma_h)\frac{d \phi_i^k}{dx}dx\\+Q(u_p^{k},u_0^{k+1},R^k)\phi_i^k(x_k)-Q(u_p^{k-1},u_0^{k},R^{k-1})\phi_i^k(x_{k-1})=0
\end{multline*}
\begin{equation*}u_0^{k+1}\phi_i^k(x_k)-u_0^k\phi_i^k(x_{k-1})-\int_{x_{k-1}}^{x_k}u_h\frac{d \phi_i^k}{dx}dx=\int_{x_{k-1}}^{x_k} \sigma_h\phi_i^k dx.\end{equation*}
Then we substitute $u_h$ and $\sigma_h$:
\begin{multline*}
    \int_{x_{k-1}}^{x_k} \partial_t \left(\sum_{j=0}^p u_j^k\phi_j^k\right)\phi_i^kdx-\int_{x_{k-1}}^{x_k}q_F(u_h,\sigma_h)\frac{d \phi_i^k}{dx}dx\\+Q(u_p^{k},u_0^{k+1},R^k)\phi_i^k(x_k)-Q(u_p^{k-1},u_0^{k},R^{k-1})\phi_i^k(x_{k-1})=0
\end{multline*}
\begin{equation*}u_0^{k+1}\phi_i^k(x_k)-u_0^k\phi_i^k(x_{k-1})-\int_{x_{k-1}}^{x_k}\left(\sum_{j=0}^p u_j^k\phi_j^k\right)\frac{d \phi_i^k}{dx}dx=\int_{x_{k-1}}^{x_k} \left(\sum_{j=0}^p \sigma_j^k\phi_j^k\right)\phi_i^k dx.\end{equation*}
We can bring out the basis coefficients, because they are constants, $\forall k,i$:
\begin{multline*}
\sum_{j=0}^p\partial_t u_j^k\int_{x_{k-1}}^{x_k} \phi_i^k\phi_j^kdx-\int_{x_{k-1}}^{x_k}q_F(u_h,\sigma_h)\frac{d \phi_i^k}{dx}dx\\+Q(u_p^{k},u_0^{k+1},R^k)\phi_i^k(x_k)-Q(u_p^{k-1},u_0^{k},R^{k-1})\phi_i^k(x_{k-1})=0
\end{multline*}
\begin{equation*}u_0^{k+1}\phi_i^k(x_k)-u_0^k\phi_i^k(x_{k-1})-\sum_{j=0}^p u_j^k\int_{x_{k-1}}^{x_k} \phi_j^k\frac{d \phi_i^k}{dx}dx=\sum_{j=0}^p\sigma_j^k\int_{x_{k-1}}^{x_k}\phi_i^k\phi_j^k dx.\end{equation*}

We now use the properties of the basis functions to simplify the boundary terms. Specifically, $\phi_i^k(x_{k-1})=1,\; i = 0$, otherwise $\phi_i^k(x_{k})=0$, and $\phi_i^k(x_{k})=1,\;i = p$, otherwise $\phi_i^k(x_{k-1})=0$. For the first equation:
\begin{align*}
    &\sum_{j=0}^p\partial_t u_j^k\int_{x_{k-1}}^{x_k} \phi_0^k\phi_j^kdx-\int_{x_{k-1}}^{x_k}q_F(u_h,\sigma_h)\frac{d \phi_0^k}{dx}dx-Q(u_p^{k-1},u_0^{k},R^{k-1})=0\\
    &\sum_{j=0}^p\partial_t u_j^k\int_{x_{k-1}}^{x_k} \phi_i^k\phi_j^kdx-\int_{x_{k-1}}^{x_k}q_F(u_h,\sigma_h)\frac{d \phi_i^k}{dx}dx=0,\quad 0<i<p\\
    &\sum_{j=0}^p\partial_t u_j^k\int_{x_{k-1}}^{x_k} \phi_p^k\phi_j^kdx-\int_{x_{k-1}}^{x_k}q_F(u_h,\sigma_h)\frac{d \phi_p^k}{dx}dx+Q(u_p^{k},u_0^{k+1},R^k)=0.
\end{align*}
For the second equation:
\begin{align*}
-u_0^k-\sum_{j=0}^p u_j^k\int_{x_{k-1}}^{x_k} \phi_j^k\frac{d \phi_0^k}{dx}dx&=\sum_{j=0}^p\sigma_j^k\int_{x_{k-1}}^{x_k}\phi_0^k\phi_j^k dx\\
-\sum_{j=0}^p u_j^k\int_{x_{k-1}}^{x_k} \phi_j^k\frac{d \phi_i^k}{dx}dx&=\sum_{j=0}^p\sigma_j^k\int_{x_{k-1}}^{x_k}\phi_i^k\phi_j^k dx,\quad 0<i<p\\
u_0^{k+1}-\sum_{j=0}^p u_j^k\int_{x_{k-1}}^{x_k} \phi_j^k\frac{d \phi_p^k}{dx}dx&=\sum_{j=0}^p\sigma_j^k\int_{x_{k-1}}^{x_k}\phi_p^k\phi_j^k dx.
\end{align*}
Turning this into a linear system, we obtain:
\begin{align*}
    &\mathbf{M}^k\mathbf{\sigma}^k = -\mathbf{C}^k\mathbf{u}^k+\mathbf{S}^k_1(u_h)\\
    &\mathbf{M}^k\dot{\mathbf{u}}^k - \mathbf{K}^k(u_h,\sigma_h)+\mathbf{S}^k_2(u_h)= 0
\end{align*}
with $(p+1)\times(p+1)$ dimension matrices $\mathbf{M}^k$ and $\mathbf{C}^k$, $(p+1)$ dimension vector-valued function $\mathbf{K}^k(u_h,\sigma_h)$ defined below, and vectors $\dot{\mathbf{u}}^k,\mathbf{u}^k,\mathbf{\sigma}^k$ that represent the real coefficients on the basis functions at the basis points. 
\begin{align*}
    M^k_{ij} &= \int_{x_{k-1}}^{x_k} \phi_i^k\phi_j^kdx\\
    C^k_{ij} &= \int_{x_{k-1}}^{x_k} \frac{d \phi_i^k}{dx}\phi_j^kdx\\
    K^k_i(u_h,\sigma_h) &= \int_{x_{k-1}}^{x_k}q_F(u_h,\sigma_h)\frac{d \phi_i^k}{dx}dx\\
    \mathbf{S}^k_1(u_h) &= \left( -u_0^k, 0, \cdots, 0, u_0^{k+1} \right)^\top\\
     \mathbf{S}^k_2(u_h) &= \left(-Q(u_p^{k-1},u_0^{k},R^{k-1}), 0, \cdots, 0, Q(u_p^{k},u_0^{k+1},R^k)\right)^\top.
\end{align*}
Since the basis functions are chosen such that they are equal to 1 at the basis points, we can view $\dot{\mathbf{u}}^k,\mathbf{u}^k,\mathbf{\sigma}^k$ as the value of the functions $\dot{u}_h,u_h,\sigma_h$ at the basis points.

Calculation of $\mathbf{K}^k$, requires the use of an integral approximation technique. For this paper, we use Gauss-Legendre quadrature with the number of points, $N_G\geq(p+1)/2$:
\begin{align*}
    K^k_i(u_h,\sigma_h) &= \int_{x_{k-1}}^{x_k}q_F(u_h,\sigma_h)\frac{d \phi_i^k}{dx}dx \\
    &= \sum_{g=1}^{N_G}w^k_gD(u_h(x^k_g))U(R(x^k_g))\frac{d \phi_i^k}{dx}(x^k_g)
\end{align*}
where $x^k_g$ and $w^k_g$ are the quadrature points and associated weights. The approximation is exact for polynomials of degree $p$, which is the case for each interval $k$. $\mathbf{S}^k_2$ and $\mathbf{K}^k$ requires approximating the non-local convolution terms for $R^k$ and $R(x_g)$:
\begin{align*}
    R^k &= \int_{x_k}^{x_k+\gamma}\mathcal K_\gamma\left(y-x_k\right)\left(u_h(y)+\kappa\Psi(\sigma_h(y))\right)dy\\
    R(x_g) &= \int_{x_g}^{x_g+\gamma}\mathcal K_\gamma\left(y-x_g\right)\left(u_h(y)+\kappa\Psi(\sigma_h(y))\right)dy.
\end{align*}

In the case of $R(x_g)$, this needs to be done for each Gauss-Legendre point. Furthermore, note that $N_G$ points are no longer sufficient for evaluating the convolution, as the integration is computed over a piecewise polynomial across multiple intervals. Instead, the integral over $[x,x+\gamma]$ is computed using piecewise integrals, split at each domain partition. We refer to \cite{chalons2018high} for the detailed formulation of the calculation of the non-local convolutions. This introduces $\mathcal O (n^2)$ complexity and computational cost (see Sections \ref{sec:complex} and \ref{sec:time}). 

At each time step, we solve for each $\dot{\mathbf u}^k$ by first solving for $\mathbf \sigma^k$ using a linear solver, such as the backslash operator, which we then use to solve for each $\dot{\mathbf u}^k$, giving $\dot{u}_h$, constructed similarly to \eqref{eq:approx_sol}:
\begin{align}
\label{eq:lin_sys_1}
    \sigma^k &= \mathbf{M}^k \setminus \left(-\mathbf C^k\mathbf u^k+\mathbf S^k_1(u_h)\right)\\
    \label{eq:lin_sys_2}
    \dot{\mathbf u}^k &= \mathbf M^k \setminus \left(\mathbf K^k(u_h,\sigma_h)-\mathbf S^k_2(u_h)\right).
\end{align}
 Finally, the discrete advancement in time uses the third-order Runge-Kutta method.
\begin{align*}
    u^{(1)} &= u^t_h + \Delta t \dot{u}^t_h\\
    u^{(2)} &= \frac{3}{4}u^t_h + \frac{1}{4}\left(u^{(1)}+\Delta t\dot{u}^{(1)}\right)\\
    u^{t+1}_h &= \frac{1}{3}u^t_h + \frac{2}{3}\left(u^{(2)}+\Delta t\dot{u}^{(2)}\right).
\end{align*}

To satisfy the CFL condition \cite{cockburn2001runge}, time step $\Delta t$ satisfies conditions: \begin{equation*}
\frac{\Delta t}{\Delta x} = \beta\frac{1}{(2p+1)\max\left|\pp_\rho(D(\rho)U(\hat \rho))\right|},
\end{equation*}
where $\beta \leq 1$ is the CFL number. Note that for our choice of $\hat \rho$, the maximum will coincide with $\max\left|\pp_\rho(D(\rho)U(\rho))\right|$.
\subsection{Generalized slope limiter}

Finally, we use a generalized slope limiter (GSL) to control oscillations, which can form when there are discontinuities in the solution (see \cite{cockburn2001runge}). This limiter is applied to all intermediate steps of the time discretization when the solution has large discontinuties. We consider the following GSL proposed in \cite{cockburn2001runge}. We define:
\begin{equation*}
    \Delta_{+}\bar u^k = \bar u^{k+1} - \bar u^k, \quad \Delta_{-}\bar u^k = \bar u^k - \bar u^{k-1}
\end{equation*}
where $\bar u^k$ is the average value of $u_h$ on interval $k$. Again, this can be computed exactly with Gauss-Legendre quadrature. We set
\begin{align*}
    u_-^{k} &= \bar u^k + \text{minmod}_2\left(u^k_p - \bar u^k,\Delta_{+}\bar u^k,\Delta_{-}\bar u^k\right)\\
    u_+^{k} &= \bar u^k - \text{minmod}_2\left(\bar u^k - u^k_0 ,\Delta_{+}\bar u^k,\Delta_{-}\bar u^k\right)
\end{align*}
where $\text{minmod}_2$ is given by the TVB modified minmod function,
\begin{align*}
    \text{minmod}_2(a,b,c)&=
    \begin{cases}
        a & \abs{a}\leq M(\Delta x)^2 \\
        \text{minmod}_1(a,b,c) & \text{otherwise}
    \end{cases}\\
    \text{minmod}_1(a,b,c)&=
    \begin{cases}
        \text{sign}(a)\cdot\min(\abs{a},\abs{b},\abs{c}) & \text{sign}\left(a\right)=\text{sign}\left(b\right)=\text{sign}\left(c\right) \\
        0 & \text{otherwise}
    \end{cases}
\end{align*}
with $M>0$ a constant. An $M$ that is too small can introduce numerical diffusion into the scheme, while a value too large may not control the oscillations (see \cite{cockburn2001runge} for discussion of $M$). Finally, the GSL replaces $u^k_h$ on each cell with the following:
\begin{equation*}
    \begin{cases}
    u^k_h & u^k_+ = u^k_0 \text{ and } u^k_- = u^k_p\\
    \bar u^k + \left(\frac{(x-x^k)}{\Delta x/2} + 1\right)\text{minmod}_2\left(u^k_1 ,\Delta_{+}\bar u^k,\Delta_{-}\bar u^k\right) & \text{otherwise}
    \end{cases}.
\end{equation*}
This GSL is used after each inner step of the time advancement. 

\subsection{Complexity of non-local LDG}
\label{sec:complex} The standard 1D DG scheme is $\mathcal{O}(n)$, as constructing each $K^k$ and $S^k_2$ takes a constant number of points dependent on $p$. Solving the linear system \eqref{eq:lin_sys_2} is also $\mathcal{O}(n)$, as the dimension is based on $p$, i.e. $\mathbf{M}^k$ is $(p+1)\times(p+1)$. Incorporating diffusion involves constructing $S^k_1$ and solving an additional equation \eqref{eq:lin_sys_1} for $\sigma_h$, which maintains linear complexity. It is the computations of the non-local convolutions which increase the order of complexity.

Under the same assumptions as \cite{chalons2018high} where $\gamma n/l\in\mathbb{Z}$ with $l$ defined as the domain length, the convolutions, $R^k$, can each be computed exactly with $N_G\times \gamma n/l$ points. Similarly, each $R(x_g)$ can be computed exactly with $N_G\times(\gamma n/l + 1)$ points. For each time step, there are $n$ number of $R^k$ which must be computed, while there are $N_G\times n$ number of $R(x_g)$. This makes the computational complexity of computing $\mathbf{K}^k$, $\mathbf{S}^k_2$ for all $k$:
\begin{align*}
    \mathbf{K}^k&: (N_G\times n)\times(N_G\times(\gamma n/l + 1)) = \mathcal{O}(n^2)\\
    \mathbf{S}^k_2&: n\times(N_G\times\gamma n/l) = \mathcal{O}(n^2).
\end{align*}
In Section \ref{sec:time}, we will show that this quadratic order is expressed in the computation time.

\section{Results}
\label{sec:results}
For the following results, we simulate the proposed model \eqref{E:main}\eqref{E:fictitious_density} on $\Omega = [0,1]$. The numerical method was coded in Julia and run on a laptop with an Intel Core i7-8650U CPU with 16 GB of available RAM. 
\subsection{Model Properties}
\label{sec:model_prop}
To evaluate the proposed model against the conventional non-local conservation law model and conventional conservation law model, we simulate the proposed model \eqref{E:main}\eqref{E:fictitious_density} with four different sets of parameters, $\{\gamma=0,\kappa=0\},\{\gamma=0.1,\kappa=0\},\{\gamma=0.1,\kappa=0.25\},\text{and } \{\gamma=0.1,\kappa=0.5\}$. The last two models are used to evaluate the effect of the diffusion intensity. For each of these models, we simulate using polynomial degree $p=1$ with $\Delta x = 1/320$ until final time $T=1$. The CFL number is set $\beta=0.2$ for these simulations and $M=35$. The boundary conditions are set to be frozen, i.e., ghost cells are defined such that $\rho(x>1,t) = \rho_\circ(1)$, $\rho(x<0,t) = \rho_\circ(0)$. 

\subsubsection{Rarefaction}
\label{sec:rarefaction}
A rarefaction scenario is simulated using the following initial condition: 
\begin{equation*}
\rho_\circ (x) =
    \begin{cases}
        0.45 & x < 0.5 \\
        0.20 & \text{otherwise}.
    \end{cases}
\end{equation*}

\begin{figure}[htbp]
  \label{fig:sol_045-020}
  \centering
  \includegraphics[width=0.7\textwidth]{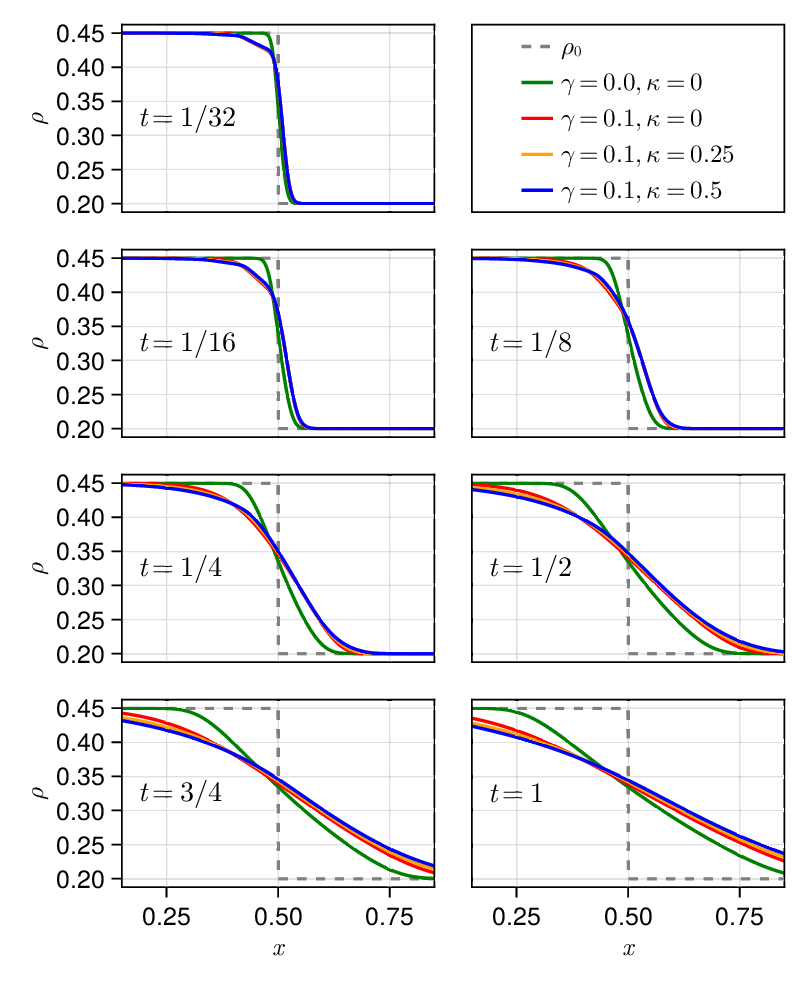}
  \caption{Rarefaction: Comparison between local, non-local, and non-local diffusive model at several timesteps. The solution was simulated with polynomial degree $p=1$ and $\Delta x = 1/320$. }
\end{figure}

In Figure \ref{fig:sol_045-020}, we see that intensifying the diffusion coefficient $\kappa$ increases the amount of flow traveling rightward. The curves associated with $\kappa>0$ have values of $\rho$ closer to $\arg\max q_F(\rho) \approx 0.33 $ across the entire solution. This is expected as $\hat\rho$ decreases due to the negative value of $\pp_x \rho$, which relatively increases the speed at which it travels rightward.

\subsubsection{Shockwaves}
\label{sec:shockwaves}
\begin{figure}[htbp]
  \label{fig:sol_015-045}
  \centering
\includegraphics[width=0.7\textwidth]{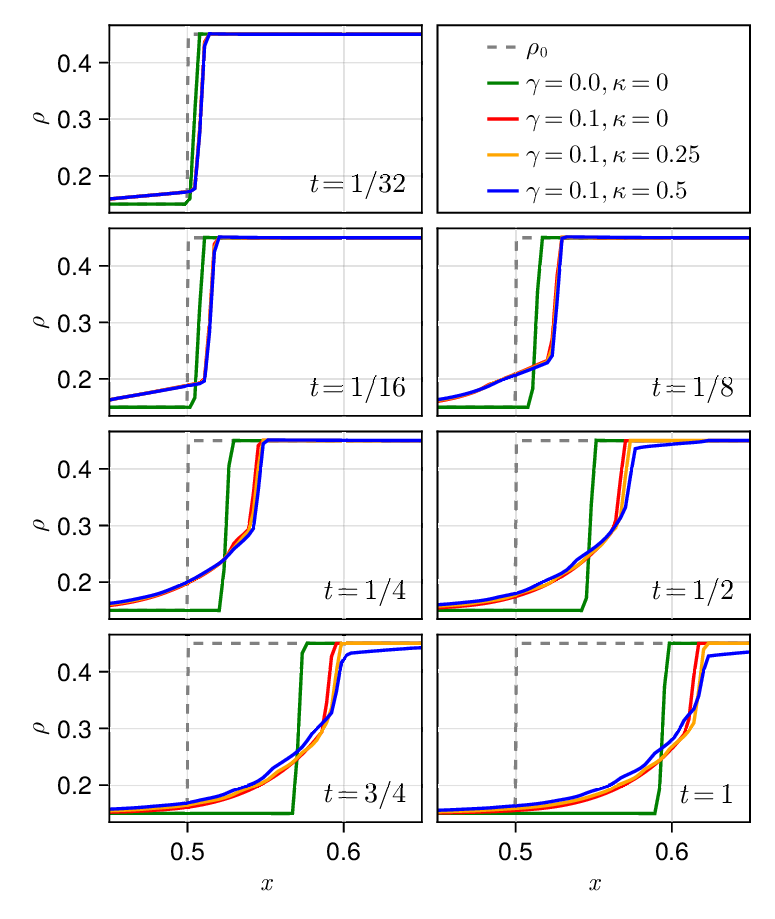}
  \caption{Forwards Shockwave: Comparison between local, non-local, and non-local diffusive model at several timesteps. The solution was simulated with polynomial degree $p=1$ and $\Delta x = 1/320$.}
\end{figure}

\begin{figure}[htbp]
  \label{fig:sol_035-065}
  \centering
  \includegraphics[width=0.7\textwidth]{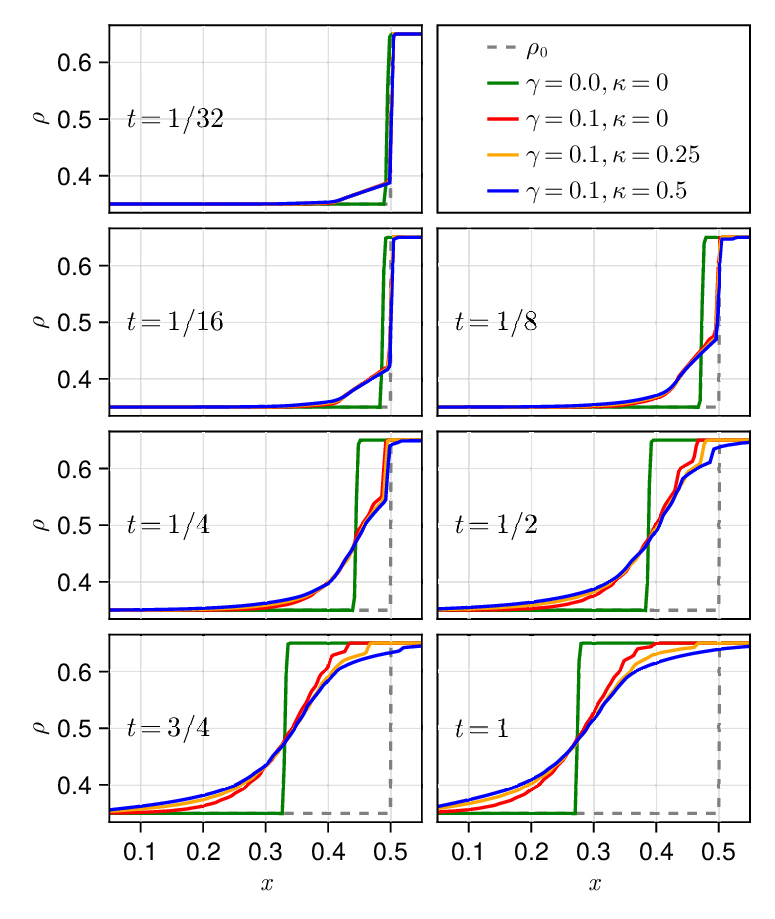}
  \caption{Backwards Shockwave: Comparison between local, non-local, and non-local diffusive model at several timesteps. The solution was simulated with polynomial degree $p=1$ and $\Delta x = 1/320$.}
\end{figure}

Two shockwave scenarios are considered, forwards shockwave and backwards shockwave. The forwards uses the following initial condition:
\begin{equation*}
\rho_\circ (x) =
    \begin{cases}
        0.15 & x < 0.5 \\
        0.45 & \text{otherwise}.
    \end{cases}
\end{equation*}
The backwards shockwave uses the following initial condition:
\begin{equation*}
\rho_\circ (x) =
    \begin{cases}
        0.35 & x < 0.5 \\
        0.65 & \text{otherwise}.
    \end{cases}
\end{equation*}
In Figure \ref{fig:sol_015-045} and \ref{fig:sol_035-065}, we see that the intensity of the diffusion coefficient contributes to smoothing out the solution compared to that of the conventional non-local model. In the shockwave cases, this can be attributed to $\hat\rho$ being greater than $\rho$, thus the solution of the diffusive models lagging behind the solution of the non-local model (seen on the left areas of the solution). Additionally, density of the diffusive models will be less than density of the non-diffusive non-local model on the right because density has not traveled rightward with as much velocity.

\begin{figure}[htbp]
  \label{fig:0-1}
  \centering
  \includegraphics[width=0.7\textwidth]{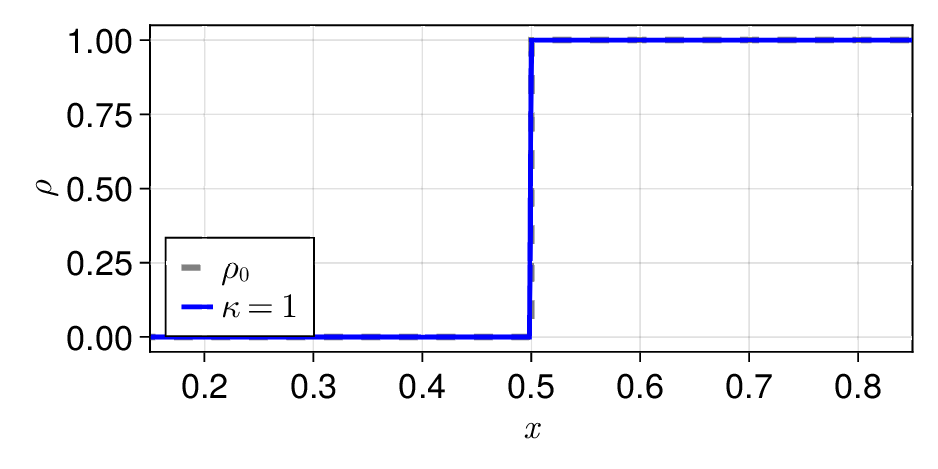}
  \caption{Absence of backwards velocity. The solution was simulated with polynomial degree $p=1$ and $\Delta x = 1/320$ with $\gamma=0.1,\kappa=1$.}
\end{figure}

Lastly, we check the special case of 
\begin{equation*}
\rho_\circ (x) =
    \begin{cases}
        0 & x < 0.5 \\
        1 & \text{otherwise}
    \end{cases}
\end{equation*}
to verify that under our construction, there is no presence of negative velocity mentioned in \cite{daganzo1995requiem}. Unlike what would happen with a linear diffusion model, our model does not change with time, which preserves one of the physical constraints of traffic flow (Figure \ref{fig:0-1}).

\subsection{Computational Results}
\label{sec:comp_results}
We simulate the solution with initial condition:
\begin{equation*}
    \rho_\circ(x) = 0.5 + 0.4\sin(2\pi (x+0.5))
    \label{eq:init_num}
\end{equation*}
until $T=0.1$. We use periodic boundary conditions and remove the generalized slope limiter, as the solution will not have discontinuities. For each test, we simulate using polynomial degrees $p=1,2,3$ and vary $n = l/\Delta x = 20, 40, 80, 160, 320$. The CFL number is set $\beta = 0.1$ for these simulations. We take the average computation time from 20 simulations for the time readings. 

\subsubsection{Error Convergence}
\label{sec:converge}

For the convergence tests, we use the solution with $p=4$ on a fine mesh where $n = 640$ as the reference solution, for a total of $(p+1)\times n = 3200$ degrees of freedom. We use the $L^2$-error:
\begin{equation*}
    L^2(\Delta x) = \sqrt{\int_{\Omega} \left(u_h - u_{\text{ref}}\right)^2dx}.
\end{equation*}
Note that with the piecewise-polynomial properties of $u_h$, and in this case $u_{\text{ref}}$, this can be computed exactly with Gauss-Legendre quadrature. 

\begin{figure}[htbp]
  \centering
  \includegraphics[width=0.7\textwidth]{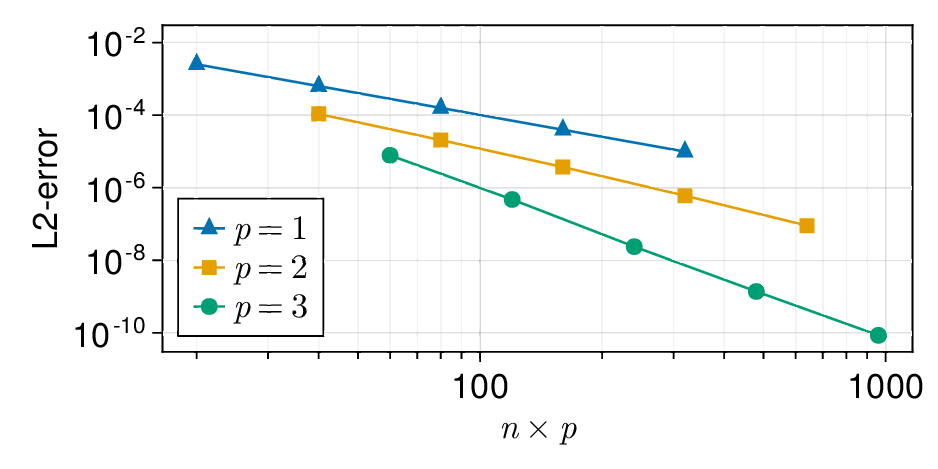}
  \caption{Degrees of freedom versus $L^2$-errors. Numerical solution for
proposed model with $\Delta x = 1/20, 1/40, 1/80, 1/160, 1/320$ at $T = 0.1$, for polynomial degrees $p=1,2,3$. The reference solution is taken as the solution for $p=4$ with $\Delta x = 1/640$}
  \label{fig:error_converge}
\end{figure}

From Figure \ref{fig:error_converge}, we can observe that higher-order approximations are computationally more efficient. The fewest tested degrees of freedom for $p=3$ has a lower error than and is faster than the highest tested degrees of freedom for $p=1$. 

We can estimate the convergence rates (Table \ref{tab:compdata}) by calculating the slopes of the lines in Figure \ref{fig:error_converge} with
\begin{equation*}
    m = \frac{\log\left[L^2(1/20)/L^2(1/320)\right]}{\log\left[(1/20)/(1/320)\right]}.
\end{equation*}
We find that the rates are similar to those found in \cite{chalons2018high} for their scheme which excludes diffusion, which is expected. 

\begin{table}[htbp]
{
    \begin{center}
        \caption{$L^2$-errors, computation time, and rates of convergence for polynomial degrees $p=1,2,3$ without slope limiter at $T=0.1$. Reference solution is taken to be solution for $p=4$, $n = 640$.}
    \label{tab:compdata}
    \begin{tabular}{|c|c|r|c|r|c|r|}
    \hline
         & \multicolumn{2}{c|}{$p=1$}  & \multicolumn{2}{c|}{$p=2$}  & \multicolumn{2}{c|}{$p=3$}  \\
         \hline
         & \multicolumn{2}{c|}{$m=1.9987$} & \multicolumn{2}{c|}{$m=2.5582$} & \multicolumn{2}{c|}{$m=4.1196$}  \\
         \hline
        $n$ & $L^2$ & CPU [s] &  $L^2$ & CPU [s]  & $L^2$ & CPU [s]\\
        \hline
        20 &  2.52e$-03$ & 0.131 & 1.08e$-04$ & 0.166 & 7.79e$-06$ & 0.324\\
        40 &  6.31e$-04$ & 0.342 & 2.05e$-05$ & 0.611 & 4.73e$-07$ & 0.725\\
        80 &  1.58e$-04$ & 1.141 & 3.71e$-06$ & 1.728 & 2.37e$-08$ & 1.837\\
        160 & 3.95e$-05$ & 3.815 & 6.03e$-07$ & 4.767 & 1.38e$-09$ & 6.341\\
        320 & 9.88e$-06$ & 15.122  & 8.99e$-08$ & 20.925 & 8.53e$-11$ & 29.301\\
        \hline
    \end{tabular}
    \end{center}
    }
\end{table}

\subsubsection{Computation Time}
\label{sec:time}

\begin{table}[htbp]
{
    \begin{center}
        \caption{Percentage of computation time to compute standard conservation law components, diffusion components, and non-local components for polynomial degree $p=3$ without slope limiter at $T=0.1$.}
    \label{tab:timebreak}
    \begin{tabular}{|c|r|r|r|}
    \hline
        & Standard & Diffusion &  Non-local \\
        $n$ & \% & \% &  \% \\
        \hline
        20 & 32.9 & 18.9 & 48.2\\
        40 & 28.1 & 12.7 & 59.2\\
        80 & 21.3 & 9.1 &  69.6\\
        160 & 14.4 & 7.9 & 77.5\\
        320 &  9.8 & 7.0 & 83.2\\
        \hline
    \end{tabular}
    \end{center}
    }
\end{table}

\begin{figure}[htbp]
  \centering
  \includegraphics[width=0.7\textwidth]{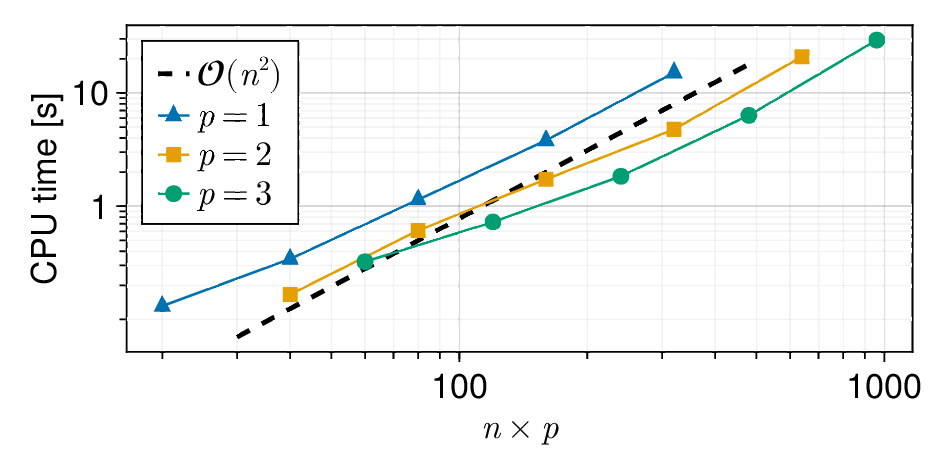}
  \caption{Degrees of freedom versus CPU time. Numerical solution for
proposed model with $\Delta x = 1/20, 1/40, 1/80, 1/160, 1/320$ at $T = 0.1$, for polynomial degrees $p=1,2,3$}
  \label{fig:time_scaling}
\end{figure}

 While this implementation has not been optimized for speed, as stated in Section \ref{sec:complex}, calculating the convolution terms is expected to increase the complexity to the order of $\mathcal O (n^2)$ due to the additional number of points that must be sampled for the convolutions. In Figure \ref{fig:time_scaling}, we can see that the computation time increases quadratically as $n$ increases for all polynomial degrees.

\begin{figure}[htbp]
  \centering
  \includegraphics[width=0.85\textwidth]{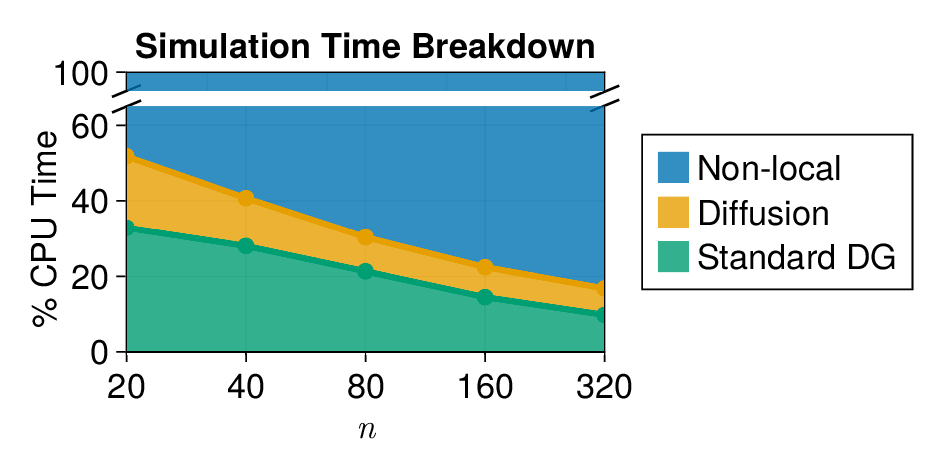}
  \caption{Cost of diffusion and non-local calculations. Partitions versus percentage of CPU Time. Benchmarks of computation time for numerical solution for
proposed model with $\Delta x = 1/20, 1/40, 1/80, 1/160, 1/320$, for polynomial degree $p=3$}
  \label{fig:diffusion_cost}
\end{figure}

To understand the computational cost of simulating a model with diffusive components and non-local convolution, we examine the time difference between simulating the proposed model with $\{\gamma > 0,\kappa > 0\}$, with $\{\gamma > 0, \kappa=0\}$, and with $\{\gamma=\kappa=0\}$. Note that the second and third simulations reduce to the non-local conservation law (presented in \cite{chalons2018high}) and standard DG for conservation laws. In Figure \ref{fig:diffusion_cost}, we see diffusion increases the CPU time consistently by about 50\% from that of solving the standard conservation law. This is expected as adding diffusion to DG amounts to solving an extra equation of the same complexity for $\sigma_h$, before solving for $u_h$, but does not require computing any integrals (i.e. only $\mathbf S_1^k$ is updated with each time step). With our implementation, we find that over 45\% of CPU time is used computing the convolutions used for non-local models, which increases above 80\% with the number of partitions, supporting the analysis in Section \ref{sec:complex}. We would expect that as $n\rightarrow\infty$, nearly the entire portion of time will be used to compute the non-local components. 

For 2D non-local problems, the complexity would be increased from quadratic to quartic. Additionally, the authors of \cite{gottlich2015discontinuous} observe that the reprocessing requirements (i.e. calculation and storage of quadrature coefficients) of DG on a non-uniform 2D mesh can exceed memory capacity for large polynomial degree and small mesh sizes. 

\begin{figure}[htbp]
  \centering
  \includegraphics[width=0.7\textwidth]{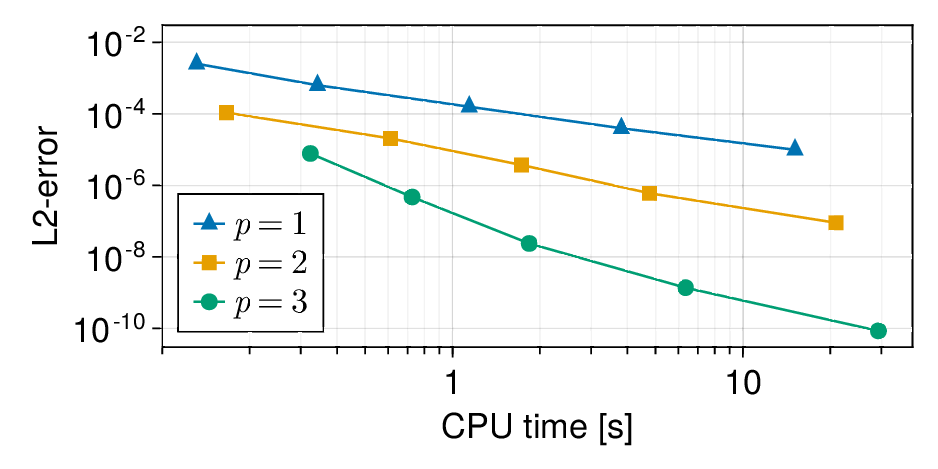}
  \caption{CPU time versus $L^2$-errors. Numerical solution for
proposed model with $\Delta x = 1/20, 1/40, 1/80, 1/160, 1/320$ at $T = 0.1$, for polynomial degrees $p=1,2,3$. The reference solution is taken as the solution for $p=4$ with $\Delta x = 1/640$}
  \label{fig:time_converge}
\end{figure}

\section{Conclusions}
\label{sec:conclusions}
In this paper, we have proposed a non-local diffusive model for traffic flow and show under what conditions it is accurately defined. A novel numerical scheme used to accurately and efficiently simulate the model was presented. The results in this paper show that the usage of higher-order schemes, like the proposed, are essential for computing accurate solutions to non-local conservation laws. Increasing the number of partitions in a 1D mesh increases the computation time at a quadratic rate, owing to the computation of the non-local convolution. However, high levels of accuracy can still be achieved in reasonable time using high-order approximations, which do not require as fine of a mesh. 

\section*{Acknowledgments}
Dawson Do is partially supported by the Dwight D. Eisenhower Transportation Fellowship Program.

\bibliographystyle{plain}
\bibliography{references}

\begin{thebibliography}{10}

\bibitem{aw2000resurrection}
AATM Aw and Michel Rascle.
\newblock Resurrection of" second order" models of traffic flow.
\newblock {\em SIAM journal on applied mathematics}, 60(3):916--938, 2000.

\bibitem{bayen2022modeling}
Alexandre Bayen, Jan Friedrich, Alexander Keimer, Lukas Pflug, and Tanya Veeravalli.
\newblock Modeling multilane traffic with moving obstacles by nonlocal balance laws.
\newblock {\em SIAM Journal on Applied Dynamical Systems}, 21(2):1495--1538, 2022.

\bibitem{betancourt2011nonlocal}
Fernando Betancourt, Raimund B{\"u}rger, Kenneth~H Karlsen, and Elmer~M Tory.
\newblock On nonlocal conservation laws modelling sedimentation.
\newblock {\em Nonlinearity}, 24(3):855, 2011.

\bibitem{blandin2016well}
Sebastien Blandin and Paola Goatin.
\newblock Well-posedness of a conservation law with non-local flux arising in traffic flow modeling.
\newblock {\em Numerische Mathematik}, 132(2):217--241, 2016.

\bibitem{bonzani2000hydrodynamic}
I~Bonzani.
\newblock Hydrodynamic models of traffic flow: Drivers' behaviour and nonlinear diffusion.
\newblock {\em Mathematical and computer modelling}, 31(6-7):1--8, 2000.

\bibitem{campos2021saturated}
Juan Campos, Andrea Corli, and Luisa Malaguti.
\newblock Saturated fronts in crowds dynamics.
\newblock {\em Advanced Nonlinear Studies}, 21(2):303--326, 2021.

\bibitem{chalons2018high}
Christophe Chalons, Paola Goatin, and Luis~M Villada.
\newblock High-order numerical schemes for one-dimensional nonlocal conservation laws.
\newblock {\em SIAM Journal on Scientific Computing}, 40(1):A288--A305, 2018.

\bibitem{cockburn2012discontinuous}
Bernardo Cockburn, George~E Karniadakis, and Chi-Wang Shu.
\newblock {\em Discontinuous Galerkin methods: theory, computation and applications}, volume~11.
\newblock Springer Science \& Business Media, 2012.

\bibitem{cockburn1989tvb2}
Bernardo Cockburn, San-Yih Lin, and Chi-Wang Shu.
\newblock Tvb runge-kutta local projection discontinuous galerkin finite element method for conservation laws iii: one-dimensional systems.
\newblock {\em Journal of computational Physics}, 84(1):90--113, 1989.

\bibitem{cockburn1989tvb}
Bernardo Cockburn and Chi-Wang Shu.
\newblock Tvb runge-kutta local projection discontinuous galerkin finite element method for conservation laws. ii. general framework.
\newblock {\em Mathematics of computation}, 52(186):411--435, 1989.

\bibitem{cockburn1998local}
Bernardo Cockburn and Chi-Wang Shu.
\newblock The local discontinuous galerkin method for time-dependent convection-diffusion systems.
\newblock {\em SIAM Journal on Numerical Analysis}, 35(6):2440--2463, 1998.

\bibitem{cockburn2001runge}
Bernardo Cockburn and Chi-Wang Shu.
\newblock Runge--kutta discontinuous galerkin methods for convection-dominated problems.
\newblock {\em Journal of scientific computing}, 16:173--261, 2001.

\bibitem{colombo2012class}
Rinaldo~M Colombo, Mauro Garavello, and Magali L{\'e}cureux-Mercier.
\newblock A class of nonlocal models for pedestrian traffic.
\newblock {\em Mathematical Models and Methods in Applied Sciences}, 22(04):1150023, 2012.

\bibitem{colombo2015nonlocal}
Rinaldo~M Colombo and Francesca Marcellini.
\newblock Nonlocal systems of balance laws in several space dimensions with applications to laser technology.
\newblock {\em Journal of Differential Equations}, 259(11):6749--6773, 2015.

\bibitem{corli2021wavefronts}
Andrea Corli and Luisa Malaguti.
\newblock Wavefronts in traffic flows and crowds dynamics.
\newblock In {\em Anomalies in Partial Differential Equations}, pages 167--189. Springer, 2021.

\bibitem{daganzo1995requiem}
Carlos~F Daganzo.
\newblock Requiem for second-order fluid approximations of traffic flow.
\newblock {\em Transportation Research Part B: Methodological}, 29(4):277--286, 1995.

\bibitem{de1999nonlinear}
Elena De~Angelis.
\newblock Nonlinear hydrodynamic models of traffic flow modelling and mathematical problems.
\newblock {\em Mathematical and computer modelling}, 29(7):83--95, 1999.

\bibitem{goatin2016well}
Paola Goatin and Sheila Scialanga.
\newblock Well-posedness and finite volume approximations of the lwr traffic flow model with non-local velocity.
\newblock {\em Networks and Heterogeneous Media}, 11(1):107--121, 2016.

\bibitem{gottlich2015discontinuous}
Simone G{\"o}ttlich, Patrick Schindler, et~al.
\newblock Discontinuous galerkin method for material flow problems.
\newblock {\em Mathematical Problems in Engineering}, 2015, 2015.

\bibitem{gugat2016analysis}
Martin Gugat, Alexander Keimer, G{\"u}nter Leugering, and Zhiqiang Wang.
\newblock Analysis of a system of nonlocal conservation laws for multi-commodity flow on networks.
\newblock {\em Networks and Heterogeneous Media}, 10(4):749--785, 2016.

\bibitem{harten1997uniformly}
Ami Harten, Bjorn Engquist, Stanley Osher, and Sukumar~R Chakravarthy.
\newblock {\em Uniformly high order accurate essentially non-oscillatory schemes, III}.
\newblock Springer, 1997.

\bibitem{huang2022stability}
Kuang Huang and Qiang Du.
\newblock Stability of a nonlocal traffic flow model for connected vehicles.
\newblock {\em SIAM Journal on Applied Mathematics}, 82(1):221--243, 2022.

\bibitem{keimer2023nonlocal}
Alexander Keimer and Lukas Pflug.
\newblock Nonlocal balance laws--an overview over recent results.
\newblock {\em Handbook of Numerical Analysis}, 24:183--216, 2023.

\bibitem{lebacque1996godunov}
JP~Lebacque.
\newblock The godunov scheme and what it means for first order traffic flow models.
\newblock In {\em Proceedings of the 13th International Symposium on Transportation and Traffic Theory, Lyon, France, July}, volume 2426, 1996.

\bibitem{leveque2002finite}
Randall~J LeVeque.
\newblock {\em Finite volume methods for hyperbolic problems}, volume~31.
\newblock Cambridge university press, 2002.

\bibitem{lighthill1955kinematic}
Michael~James Lighthill and G~Be Whitham.
\newblock On kinematic waves i. flood movement in long rivers.
\newblock {\em Proceedings of the Royal Society of London. Series A. Mathematical and Physical Sciences}, 229(1178):281--316, 1955.

\bibitem{lighthill1955kinematic2}
Michael~James Lighthill and Gerald~Beresford Whitham.
\newblock On kinematic waves ii. a theory of traffic flow on long crowded roads.
\newblock {\em Proceedings of the Royal Society of London. Series A. Mathematical and Physical Sciences}, 229(1178):317--345, 1955.

\bibitem{matin2023Nonlinear}
Hossein Nick Zinat~Matin, Dawson Do, and Maria~Laura {Delle Monache}.
\newblock Nonlinear advection-diffusion models of traffic flow: a numerical study.
\newblock In {\em Proceedings IEEE-ITSC-The 26th International IEEE Conference on Intelligent Transportation Systems}. IEEE Society, in press.

\bibitem{pareschi2005implicit}
Lorenzo Pareschi and Giovanni Russo.
\newblock Implicit--explicit runge--kutta schemes and applications to hyperbolic systems with relaxation.
\newblock {\em Journal of Scientific computing}, 25:129--155, 2005.

\bibitem{perthame2006transport}
Beno{\^\i}t Perthame.
\newblock {\em Transport equations in biology}.
\newblock Springer Science \& Business Media, 2006.

\bibitem{pflug2020emom}
Lukas Pflug, Tobias Schikarski, Alexander Keimer, Wolfgang Peukert, and Michael Stingl.
\newblock emom: Exact method of moments—nucleation and size dependent growth of nanoparticles.
\newblock {\em Computers \& Chemical Engineering}, 136:106775, 2020.

\bibitem{richards1956shock}
Paul~I Richards.
\newblock Shock waves on the highway.
\newblock {\em Operations research}, 4(1):42--51, 1956.

\bibitem{shu1988efficient}
Chi-Wang Shu and Stanley Osher.
\newblock Efficient implementation of essentially non-oscillatory shock-capturing schemes.
\newblock {\em Journal of computational physics}, 77(2):439--471, 1988.

\bibitem{shu1989efficient}
Chi-Wang Shu and Stanley Osher.
\newblock Efficient implementation of essentially non-oscillatory shock-capturing schemes, ii.
\newblock {\em Journal of Computational Physics}, 83(1):32--78, 1989.

\bibitem{zhang2002non}
H~Michael Zhang.
\newblock A non-equilibrium traffic model devoid of gas-like behavior.
\newblock {\em Transportation Research Part B: Methodological}, 36(3):275--290, 2002.

\end{thebibliography}
\end{document}